\documentclass[letterpaper, 10 pt, conference]{ieeeconf}  

\IEEEoverridecommandlockouts                              
\usepackage{graphics} 
\usepackage{epsfig} 
\usepackage{times} 
\usepackage{amsmath} 
\usepackage{amssymb}  
\usepackage{changes}
\usepackage{booktabs}
\usepackage{makecell}
\usepackage{subfig}
\usepackage{overpic}
\usepackage{enumerate}
\usepackage{algorithm}

\newtheorem{theorem}{Theorem}
\newtheorem{proposition}{Proposition}

\newtheorem{corollary}{Corollary}
\newtheorem{assumption}{Assumption}

\newcommand{\mb}[1]{\mathbb{#1}}
\newcommand{\mc}[1]{\mathcal{#1}}
\newcommand{\mbf}[1]{\mathbf{#1}}
\newcommand{\N}{\mb{N}}

\newcommand{\R}{\mb{R}}
\newcommand{\Rn}[1]{\R^{#1}}
\newcommand{\Rnm}[2]{\R^{{#1}\times{#2}}}
\newcommand{\support}[2]{\operatorname{h}\left(#1,#2\right)}
\newcommand{\supportbb}[1]{\operatorname{h}(#1,\cdot)}
\newcommand{\pf}{\textbf{Proof. }}
\newcommand{\qed}{\hfill$\square$}

\newif\ifRedify
\Redifytrue 
 
\title{\LARGE \bf
Robustifying Model Predictive Control of\\  Uncertain Linear Systems with Chance Constraints
}

\author{Kai Wang, Kiet Tuan Hoang and S\'ebastien Gros
\thanks{$^*$ The authors gratefully acknowledge support from Norwegian University of Science and Technology. After the initial submission, K. Wang became a research associate at Loughborough University and acknowledges support from Loughborough University through Grant EP/T005734/1. }
\thanks{K. Wang, K.T. Hoang and S. Gros are with the Norwegian University of Science and Technology, Trondheim 7491, Norway (e--mail: kai.wang@ntnu.no; kiet.t.hoang@ntnu.no; sebastien.gros@ntnu.no).}	
}

\begin{document}

\maketitle
\thispagestyle{empty}
\pagestyle{empty}

\begin{abstract}
This paper proposes a model predictive controller for discrete-time linear systems with additive, possibly unbounded, stochastic disturbances and subject to chance constraints. By computing a polytopic probabilistic positively invariant set for constraint tightening with the help of the computation of the minimal robust positively invariant set, the chance constraints are guaranteed, assuming only the mean and covariance of the disturbance distribution are given. The resulting online optimization problem is a standard strictly quadratic programming, just like in conventional model predictive control with recursive feasibility and stability guarantees and is simple to implement. A numerical example is provided to illustrate the proposed method.
\end{abstract}

\section{INTRODUCTION}
\label{sec:intro}

In the last two decades, model predictive control (MPC) has emerged as one of the most prevalent optimization-based control techniques in academia and industry~\cite{Morari1999,Qin2003,Racovic:2018}. The main advantages of MPC come from its ability to simultaneously cope with constraints, optimize performance and ensure desirable control-theoretic properties~\cite{Mayne2000,rawlings:mayne:2009}. However, many of these properties deteriorate in the presence of disturbances, although nominal (certainty-equivalent) MPC exhibits inherent robustness to some extent~\cite{Grimmand2007}. Nevertheless, an MPC that explicitly considers
uncertainty is required in safety-critical applications. When the uncertainty description of the disturbance is available, robust MPC~\cite{rakovic:2015,Houska2019}, stochastic MPC~\cite{kouvaritakis2015model}, and more recent distributionally robust MPC~\cite{parys:2016} can be employed to maintain safety to certain acceptable levels. Whereas robust MPC considers the set-membership description of bounded uncertainties,  stochastic MPC and distributionally robust MPC consider uncertainties which follow specific probability distributions or sets of probability distributions, respectively. 

To alleviate the computational complexity of the exact min-max robust MPC~\cite{rawlings:mayne:2009}, tube MPC has been widely investigated as a sensible option for its intuitive simplicity, computational practicality and guaranteed control-theoretic properties. Building on the two early tube MPC methods for linear systems under bounded disturbance~\cite{mayne:langson:2001,chisci:rossiter:zappa:2001}, which were proposed almost simultaneously, numerous subsequent studies have emerged. These studies have explored tube-based robust MPC~\cite{mayne:seron:rakovic:2005,rakovic:kouvaritakis:findeisen:cannon:2011,rakovic:levine:acikmese:2016,Rakovic2013}, tube-based stochastic MPC~\cite{kouvaritakisk2010,Cannon2009,Cannon2011Tac,Lorenzen2017,Farina2013,Hewing2018}, and tube-based distributionally robust MPC~\cite{wu2022,Aolaritei2023,zhong2023}.  Interestingly, most tube-based robust MPC proposals reviewed above were generated from~\cite{mayne:langson:2001}, while almost all tube-base stochastic and distributionally robust MPC proposals have adopted method~\cite{chisci:rossiter:zappa:2001}. The key difference lies in the initialization of the optimal control problem (OCP): in \cite{mayne:langson:2001}, the OCP is solved using the current state of the nominal system, whereas in \cite{chisci:rossiter:zappa:2001}, it is solved based on the current state of the actual system. Another significant difference is in how the stage constraints are tightened: \cite{mayne:langson:2001} employs (minimal) robust positively invariant sets, while \cite{chisci:rossiter:zappa:2001} utilizes a sequence of reachable sets starting from the origin. Recursive feasibility and stability are easier to achieve in \cite{mayne:langson:2001} than in \cite{chisci:rossiter:zappa:2001}, as argued in \cite{mayne:2018}. 

However, extending \cite{mayne:langson:2001} to linear stochastic systems with possibly unbounded stochastic disturbances requires computing a probabilistic positively invariant set~\cite{kofman2012} which is non-trivial. In this regard,~\cite{kofman2012} proposes a method to compute polyhedral probabilistic positively invariant sets by using Chebyshev inequality, assuming that the system matrix is diagonalizable and that the disturbance distribution has zero-mean and diagonal covariance matrix. The paper~\cite{Hewing2018Corr} proposes to characterize ellipsoidal probabilistic positively invariant sets, assuming that the system matrix is invertible and that the disturbance distribution has zero-mean. In~\cite{Fiacchini:2021}, the ellipsoidal probabilistic positively invariant sets are investigated, assuming a correlation bound for the system matrix.

\textbf{Contributions:} In this paper, we generalize the method of~\cite{mayne:langson:2001} to linear systems under possibly unbounded stochastic disturbance, where chance constraints on the state and input are considered. The proposed MPC only requires knowledge of the mean and covariance of the stochastic disturbance without assuming zero-mean, invertible/diagonal system matrices or bounds on the mean and covariance of the disturbance as required in~\cite{kofman2012,Hewing2018Corr,Fiacchini:2021}. To enforce the chance constraints through constraint tightening, a polytopic (minimal) probabilistic positively invariant set is computed, which is minimal for a group of polytopes that are probabilistic positively invariant. The resulting online optimization problem is a simple, standard, strictly convex quadratic programming for which recursive feasibility and stability properties can be easily guaranteed.   

\textbf{Structure:} Section~\ref{sec:02} details the problem setup on the systems and constraints. Section~\ref{sec:03} discusses the probabilistic positively invariant set and its computation. Section~\ref{sec:04} presents the proposed MPC algorithm. Section~\ref{sec:05} demonstrates the proposed approach with a numerical case study, and the paper is concluded in Section~\ref{sec:06}. 

\textbf{Notation:} The sets of real numbers, non-negative integers and positive integers
are denoted by $\R$, $\N$ and $\N_+$, respectively. Given $a,b\in\N$, with $a<b$, we use the notation $\mb I_a^b$ to denote the set of non-negative integers $\{a, a+1,\ldots, b\}$. The Minkowski sum and Pontryagin difference of nonempty sets $\mc{X}$ and $\mc{Y}$ are
\begin{align*}
	\mc{X}\oplus\mc Y&:= \left\{x+y\, : \, x\in\mc{X}, \ y\in\mc{Y}\right\}\text{ and}\\
	\mc{X}\ominus\mc Y&:= \left\{x\, : \,  \forall y\in\mc{Y}, \ x+y\in\mc X\right\},
\end{align*}
respectively. The image of a nonempty set $\mc{X}$ under a matrix of compatible dimensions $M$ is given by $M\mc{X}:=\left\{Mx \, :\, x\in \mc{X}\right\}$.
Likewise, if $M$ is a square matrix, for any integer $k\in\N$, $M^k\mc{X}:=\left\{M^kx\, :\, x\in \mc{X}\right\}$.
The sets of symmetric positive semi-definite and symmetric positive definite matrices in $\Rn{n}$ are denoted by $\mathbb S^{n}_{+}$ and $\mathbb S^{n}_{++}$, respectively. Ellipsoids with shape matrix $M\in\mathbb S_+^{n}$, center $c\in\Rn{n}$ and radius $\sqrt{r}>0$ are defined as 
\begin{align*}
	\mc E(M, c, r):=&\{x\,:\, (x-c)^\intercal M^{-1}(x-c)\leq r\}\\
	=&\{c+M^{\frac{1}{2}}x\,:\, x^\intercal x\leq r\},
\end{align*}
where the matrix $M^\frac{1}{2}$ is such that $M^\frac{1}{2}M^\frac{1}{2}=M$. In addition, it holds that $\mc E(M, c, r)=c\oplus\mc E(M, 0, r)$ by definition. The support function $\supportbb{\mc{X}}$ of a nonempty, closed, convex set $\mc{X}\subseteq\R^n$ is given, for all $y\in\R^n$, by
\begin{equation*}
	\support{\mc{X}}{y}:= \sup_x\, \{y^\intercal x\ :\ x\in \mc{X}\},
\end{equation*}

\section{PROBLEM SETUP}
\label{sec:02}
This paper considers uncertain linear systems of the form
\begin{equation}\label{eq:02:Sys}
	\forall k\in \N,\quad x_{k+1} = Ax_k+Bu_k+w_k,
\end{equation}
where $x_k\in\Rn{n}$ and $u_k\in\Rn{m}$ denote the state and input, respectively, while $w_k\in\Rn{n}$ denotes a randomly distributed disturbance. These disturbances $w_0, w_1,\ldots:\Rn{n}\to \Rn{n}$ are assumed identically distributed and independent of each other and $x_k$ and $u_k$, with common probability distribution $\omega$, i.e., for all $k\in\N$, we have that $w_k\sim\omega$. Due to the disturbance $w_k$, we consider chance constraints on the state $x_k$ and input $u_k$, given by: 
\begin{align}
	\forall k\in\N_+,\quad  &\Pr(x_k\in\mc X\mid x_0) \geq 1-\epsilon_x \text{ and }\label{eq:02:PxkX}\\ 
	\forall k\in\N,\quad  & \Pr(u_k\in\mc U)\geq 1-\epsilon_u, \label{eq:02:PukU}
\end{align}
where $\epsilon_x, \epsilon_u\in(0,1)$ are modeling parameters, $x_0$ is an initial state, and $\mc X$ and $\mc U$ denote constraint sets on the state and input respectively. The set $\mc X$ ($\mc U$) is a region for which the system state $x_k$ (input $u_k$) should be confined with a probability no lower than the specified confidence level $1-\epsilon_x$ ($1-\epsilon_u$).

We make the following standing assumptions on the system dynamics \eqref{eq:02:Sys} and constraints \eqref{eq:02:PxkX}-\eqref{eq:02:PukU}.
\begin{assumption}\label{ass:01}$\ $
	\begin{enumerate}[$\,\,\ $i.]
		\item The matrix pair $(A, B)\in\Rnm{n}{n}\times\Rnm{n}{m}$ is known exactly, and it is strictly stabilizable.
		\item The set $\mc X\subseteq\Rn{n}$ is a convex polyhedron that contains the origin in its interior, and the set $\mc U\subseteq\Rn{m}$ is a convex polytope that contains the origin in its interior. In addition, their irredundant inequalities characterizations are given by: 
		\begin{align}
			\mc X &:=\left\{x\, :\, \forall i\in\mc I_{\mc X},\ f_i^\intercal x\leq 1   \right\},\\
			\mc U &:=\left\{u\, :\, \forall i\in\mc I_{\mc U},\ g_i^\intercal u\leq 1   \right\},
		\end{align}
		with the known finite index sets $\mc I_{\mc X}:=\{1, 2, \ldots, n_\mc X\}$, $\mc I_{\mc U}:=\{1, 2, \ldots, n_\mc U\}$, and the known vectors $f_i\in\Rn{n}$ and $g_i\in\Rn{n}$.
		\item The probability distribution $\omega$ is not known exactly, however, its mean vector $\mu_\omega\in\Rn{n}$ and covariance matrix $\Sigma_\omega\in\mb S^n_{++}$ are known.
	\end{enumerate}
\end{assumption}

In what follows, we introduce the linear feedback law as
\begin{equation}\label{eq:02:controldecomposition}
	\forall k\in \N,\quad	u_k=v_k+Ks_k
\end{equation}
and rewrite the state $x_k$ as 
\begin{equation}\label{eq:02:statedecomposition}
	\forall k\in \N,\quad	x_k=z_k+s_k.
\end{equation}
Here, $z_k\in\Rn{n}$ denotes the nominal (disturbance-free) component such that 
\begin{equation}\label{eq:02:nominalSys}
	\forall k\in \N,\quad z_{k+1} = A z_k +Bv_k
\end{equation}
based on the nominal control input $v_k\in\Rn{m}$. Thus, by subtracting~\eqref{eq:02:nominalSys} from~\eqref{eq:02:Sys} and in view of~\eqref{eq:02:controldecomposition}--\eqref{eq:02:statedecomposition}, $s_k\in\Rn{n}$ denotes the error component satisfying
\begin{equation}\label{eq:02:disturbanceSys}
	\forall k\in \N,\quad	s_{k+1} = (A+BK)s_k + w_k, \text{ with } w_k\sim\omega. 
\end{equation}
The matrix $K\in\Rnm{m}{n}$ is \textit{a priori} given and constructed to satisfy the following assumption.
\begin{assumption}\label{ass:02}
	The matrix $K\in\Rnm{m}{n}$ is given, and it is such
	that the matrix $A_K:=A+BK$ is strictly stable, i.e., all eigenvalues of $A_K$ are in the open unit disc.
\end{assumption}

\section{PROBABILISTIC POSITIVELY INVARIANT SET}
\label{sec:03}
Let $\mu_k\in\Rn{n}$ and $\Sigma_k\in\mb S^n_{++}$ denote the mean vector and covariance matrix of the random vector $s_k$, whose stochasticity is caused by the random disturbance sequence $\{w_i\}_{i=0}^{k-1}$ via stochastic dynamics~\eqref{eq:02:disturbanceSys}. Then, the dynamics of $\mu_k$ and $\Sigma_{k}$ are specified, with $\mu_0=s_0$ and $\Sigma_0=0$, by
\begin{align}
	\forall k\in \N,\quad \mu_{k+1}    &= A_K\mu_k + \mu_\omega\text{ and } \label{eq:03:muk}\\
	\forall k\in \N,\quad \Sigma_{k+1} &= A_K\Sigma_{k}A_K^\intercal+\Sigma_\omega.\label{eq:03:Sigmak}
\end{align}
Based on the multivariate Chebyshev inequality~\cite{chen2011,navarro2016}, we have the following proposition.
\begin{proposition}
	Suppose Assumption~\ref{ass:01}--iii holds. Consider $\mu_k$ and $\Sigma_{k}$ generated by~\eqref{eq:03:muk} and~\eqref{eq:03:Sigmak}. It follows that
	\begin{equation} \label{eq:03:Pr_skinE_k}
		\forall k\in\N_{+},\quad \Pr\left(s_k\in\mc E(\Sigma_{k}, \mu_k, n/\epsilon)\mid s_0\right)\geq 1-\epsilon,
	\end{equation}
	i.e., $\mc E(\Sigma_{k}, \mu_k, n/\epsilon)$ are confidence regions with confidence level $1-\epsilon$ of the random vectors $s_k$.
\end{proposition}
\pf It follows from~\cite[Theorem 1]{chen2011} that 
\begin{align*}
	\forall k\in\N_{+},\quad &\Pr\left((s_k-\mu_k)^\intercal\Sigma_{k}^{-1}(s_k-\mu_k)> n/\epsilon \mid s_0\right)\\
	&\leq\Pr\left((s_k-\mu_k)^\intercal\Sigma_{k}^{-1}(s_k-\mu_k)\geq n/\epsilon \mid s_0\right)\\
	&\leq\epsilon,
\end{align*}
such that
\begin{align*}
	\forall k\in\N_{+},\   &\Pr\left((s_k-\mu_k)^\intercal\Sigma_{k}^{-1}(s_k-\mu_k)\leq n/\epsilon \mid s_0\right)\\
	&=1-\Pr\left((s_k-\mu_k)^\intercal\Sigma_{k}^{-1}(s_k-\mu_k)> n/\epsilon \mid s_0\right)\\
	&\geq 1-\epsilon.
\end{align*}
Thus, it follows from the definition of ellipsoids that
\begin{equation*}
	\forall k\in\N_{+},\quad \Pr\left(s_k\in\mc E(\Sigma_{k}, \mu_k, n/\epsilon)\mid s_0\right)\geq 1-\epsilon.
\end{equation*}
\qed


\begin{corollary}
	Suppose Assumptions~\ref{ass:01} and~\ref{ass:02} hold. Consider $\mu_k$ and $\Sigma_{k}$ generated by~\eqref{eq:03:muk} and~\eqref{eq:03:Sigmak}. It follows that 
	\begin{align}
		\forall k\in\N_+,\quad &\mc E(\Sigma_{k+1}, \mu_{k+1}, n/\epsilon)\nonumber\\
		&\subseteq A_K\mc E(\Sigma_{k}, \mu_k, n/\epsilon)\oplus\mc E(\Sigma_\omega, \mu_\omega, n/\epsilon).\label{eq:03:Ek}
	\end{align}
\end{corollary}
\pf By the property of ellipsoids, it holds that
\begin{equation*}
	\forall k\in\N_{+},\quad \mc E(\Sigma_{k+1}, \mu_{k+1}, n/\epsilon) = \mu_{k+1}\oplus\mc E(\Sigma_{k+1}, 0, n/\epsilon)
\end{equation*}
and 
\begin{align*}
	&A_K\mc E(\Sigma_{k}, \mu_k, n/\epsilon)\oplus\mc E(\Sigma_\omega, \mu_\omega, n/\epsilon)\\
	&=A_K\mu_k+\mu_\omega\oplus A_K\mc E(\Sigma_{k}, 0, n/\epsilon)\oplus\mc E(\Sigma_\omega, 0, n/\epsilon). 	
\end{align*}
For all $k\in\N_+$, it follows from~\cite[Property 1]{Fiacchini:2021} that 
\[
\mc E(\Sigma_{k+1}, 0, n/\epsilon) \subseteq A_K\mc E(\Sigma_{k}, 0, n/\epsilon)\oplus\mc E(\Sigma_\omega, \mu_\omega, n/\epsilon),
\] 
and it follows from~\eqref{eq:03:muk} that $$\mu_{k+1} = A_K\mu_k + \mu_\omega.$$ Thus, the Minkowski sum of both sides of the above two expressions completes the proof. \qed

At this point, we introduce the sets $\mc R_k^\epsilon$ defined by the following set-dynamics:
\begin{equation}\label{eq:03:Rk}
	\forall k\in\N,\quad \mc R^\epsilon_{k+1} = A_K \mc R^\epsilon_k\oplus\mc E(\Sigma_\omega, \mu_\omega, n/\epsilon)
\end{equation}
with the boundary condition given by $\mc R^\epsilon_0=\{s_0\}$. 

\begin{proposition}\label{prop:02}
	Suppose Assumptions~\ref{ass:01} and~\ref{ass:02} hold. Consider the sets $\mc R^\epsilon_k$ generated by~\eqref{eq:03:Rk} and $\mu_k$ and $\Sigma_{k}$ generated by~\eqref{eq:03:muk} and~\eqref{eq:03:Sigmak}.  It follows that
	\begin{equation}\label{eq:03:EkRk}
		\forall k\in\N_+,\quad \mc E(\Sigma_{k}, \mu_k, n/\epsilon)\subseteq \mc R^\epsilon_k,
	\end{equation}
	such that, for all $k\in\N_+$
	\begin{align}
		\Pr\left(s_k\in\mc R^\epsilon_k\mid s_0\right)&\geq\Pr\left(s_k\in\mc E(\Sigma_{k}, \mu_k, n/\epsilon)\mid s_0\right)\nonumber\\
		&\geq 1-\epsilon.\label{eq:03:PrRk}
	\end{align}
\end{proposition}
\pf We will show by induction that the inclusions~\eqref{eq:03:EkRk} holds. When $k=1$, we have that
\begin{align*}
	\mc E(\Sigma_{1}, \mu_1, n/\epsilon) &= A\mu_0\oplus\mc E(\Sigma_\omega, \mu_\omega, n/\epsilon)\\
	\mc R^\epsilon_1 &= As_0\oplus\mc E(\Sigma_\omega, \mu_\omega, n/\epsilon)
\end{align*}
such that $\mc E(\Sigma_{1}, \mu_1, n/\epsilon)=\mc R^\epsilon_1\subseteq\mc R^\epsilon_1$ since $\mu_0=s_0$. 
Suppose, for some $k\in\N_+$, $\mc E(\Sigma_{k}, \mu_k, n/\epsilon)\subseteq \mc R^\epsilon_k$, then we have that 
\begin{align*}
	\mc E(\Sigma_{k+1}, \mu_{k+1}, n/\epsilon)&\subseteq A_K\mc E(\Sigma_{k}, \mu_k, n/\epsilon)\oplus\mc E(\Sigma_\omega, \mu_\omega, n/\epsilon)\\
	&\subseteq A_K\mc R^\epsilon_k\oplus\mc E(\Sigma_\omega, \mu_\omega, n/\epsilon)\\
	&=\mc R^\epsilon_{k+1},
\end{align*}
in which the first inclusion follows from~\eqref{eq:03:Ek}. It follows from~\eqref{eq:03:Pr_skinE_k} and~\eqref{eq:03:EkRk} that~\eqref{eq:03:PrRk} holds. \qed

We recall from~\cite{kofman2012} that a set $\mc S\subseteq\Rn{n}$ is called probabilistic positively invariant for the system~\eqref{eq:02:disturbanceSys} with probability $1-\epsilon$, if and only if, for any $s_0\in\mc S$, $\Pr\left(s_k\in\mc S\mid s_0\right)\geq1-\epsilon$ for all $k\in\N_+$.

\begin{proposition}\label{prop:03}
	Suppose Assumptions~\ref{ass:01} and~\ref{ass:02} hold. Consider a set $\mc R^\epsilon\subseteq\Rn{n}$ satisfying, for any $s_0\in\mc R^\epsilon$,
	\begin{equation}\label{eq:03:AkRR}
		A_K\mc R^\epsilon\oplus\mc E(\Sigma_\omega, \mu_\omega, n/\epsilon)\subseteq\mc R^\epsilon
	\end{equation}
	with $\epsilon\in(0,1)$. Then, $\mc R^\epsilon$ is a probabilistic positively invariant set for the system~\eqref{eq:02:disturbanceSys} with probability $1-\epsilon$.
\end{proposition}
\pf Consider the sets $\mc R^\epsilon_k$ generated by~\eqref{eq:03:Rk}. As shown in~\eqref{eq:03:PrRk}, we have that
\[
\forall k\in\N_+,\quad \Pr\left(s_k\in\mc R^\epsilon_k\mid s_0\right)\geq 1-\epsilon.
\]  
To show that $\mc R^\epsilon$ is a probabilistic positively invariant set for the stochastic system~\eqref{eq:02:disturbanceSys} with probability $1-\epsilon$, we just need to prove that for all $s_0\in\mc R^\epsilon$ it holds that  $\mc R^\epsilon_k\subseteq\mc R^\epsilon,\, k\in\N$. When $k=0$, $\mc R^\epsilon_0=\{s_0\}\subseteq\mc R^\epsilon$ trivially holds. Suppose, for some $k\in\N$, $\mc R^\epsilon_k\subseteq\mc R^\epsilon$, then we have that
\begin{align*}
	\mc R^\epsilon_{k+1}&=A_K \mc R^\epsilon_k\oplus\mc E(\Sigma_\omega, \mu_\omega, n/\epsilon)\\
	&\subseteq A_K \mc R^\epsilon\oplus\mc E(\Sigma_\omega, \mu_\omega, n/\epsilon)\subseteq\mc R^\epsilon.
\end{align*} 
Thus, by induction, for all $s_0\in\mc R^\epsilon$, we have $\mc R^\epsilon_k\subseteq\mc R^\epsilon$ for all $k\in\N$.  \qed

In the rest of this section, we make the following standing assumption
\begin{assumption}\label{ass:03}
	Given $\epsilon\in(0,1)$, the set  $\mc E(\Sigma_\omega, \mu_\omega, n/\epsilon)$ contains the origin in its interior.
\end{assumption}

Note that Assumption~\ref{ass:03} implies that the support functions of $\mc E(\Sigma_\omega, \mu_\omega, n/\epsilon)$ is positive, i.e., for all $y\in\Rn{n}$, $\support{\mc E(\Sigma_\omega, \mu_\omega, n/\epsilon)}{y}>0$. We will discuss the methods for computing a probabilistic positively invariant set $\mc R^\epsilon$ based on the sufficient condition given by~\eqref{eq:03:AkRR}. 

Given $\epsilon\in(0,1)$ satisfying Assumption~\ref{ass:03}, we recall that a set $\mc R^\epsilon\subseteq\Rn{n}$ is robust positively invariant for the uncertain dynamics 
$\forall k\in \N,\,s_{k+1} = A_Ks_k + w_k$ with  $w_k\in \mc E(\Sigma_\omega, \mu_\omega, n/\epsilon)$ if and only if 
\begin{equation}
A_K\mc R^\epsilon\oplus\mc E(\Sigma_\omega, \mu_\omega, n/\epsilon)\subseteq\mc R^\epsilon,
\end{equation}
and it is the minimal robust positively invariant set if and only if
\begin{equation}\label{eq:03:minmalCond}
	A_K\mc R^\epsilon\oplus\mc E(\Sigma_\omega, \mu_\omega, n/\epsilon)=\mc R^\epsilon.
\end{equation}
Thus, condition~\eqref{eq:03:AkRR} is equivalent to the sufficient and necessary condition for characterizing robust positively invariant sets for the system 
$\forall k\in \N,\,s_{k+1} = A_Ks_k + w_k$ with  $w_k\in \mc E(\Sigma_\omega, \mu_\omega, n/\epsilon)$. Moreover, the minimal robust positively invariant set given by~\eqref{eq:03:minmalCond} is desirable since it will be used to tighten the constraint set in the next section. 

However, computing an exact representation of the minimal robust positively invariant set is generally impossible. In practice, outer invariant approximations of the minimal robust positively invariant set are computed~\cite{Rakovic2005}. In our setting, since $\mc E(\Sigma_\omega, \mu_\omega, n/\epsilon)$ is an ellipsoid rather than a polytope, the method~\cite{Rakovic2005} is not applicable. Instead, we utilize the approach from~\cite{Trodden2016}, which builds upon the method originally proposed in~\cite{Rakovic2013}. This approach computes a polytopic robust positively invariant set that is minimal with respect to a group of robust positively invariant sets defined by a finite number of inequalities with pre-specified normal vectors. In the sequel, we briefly summarize the method of~\cite{Trodden2016}. 
Consider the set $\mc R^\epsilon(q)$ defined as follows, with $\mc I_\mc P=\{1,...,r\}$ and $ r\in\N_+$: 
\begin{equation}
	\mc R^\epsilon(q) = \{x\in\Rn{n}:\forall i\in\mc I_\mc P, \ p_i^\intercal x \leq q_i\},
\end{equation}
where $\{p_i:\forall i\in\mc I_\mc P\}$ spans $\Rn{n}$, $q_i\geq0$ for all $i\in\mc I_\mc P$, and $q=(q_1,\ldots,q_r)^\intercal$. Then, the polytopic, minimal robust positively invariant set generated from the pre-defined normal vectors $\{p_i:\forall i\in\mc I_\mc P\}$ is given by 
\begin{equation}\label{eq:04:Rq}
	\mc R^\epsilon(q^*) = \{x\in\Rn{n}:\forall i\in\mc I_\mc P, \ p_i^\intercal x \leq q^*_i\}.
\end{equation}
Here, $q^*=d^* + c^*$, where $d^*:=(d_1^*,...,d_r^*)^\intercal$ is given by 
\begin{align}
	\forall i\in\mc I_\mc P,\ d_i^*&=\support{\mc E(\Sigma_\omega, \mu_\omega, n/\epsilon)}{p_i}\nonumber\\
	&= \mu_\omega^\intercal p_i+\sqrt{\frac{n(p_i^\intercal\Sigma_\omega p_i)}{\epsilon}},\label{eq:04:dstar}
\end{align}
and $c^*:=(c_1^*,...,c_r^*)^\intercal$ is given by the solution to the following linear programming problem:
\begin{equation}\label{eq:04:cstar}
	\begin{alignedat}{2}
		&\max_{c,\, \xi} &&\  \sum^{r}_{i=0} c_i  \\
		&\,\,\,\,\text{s.t.} &&\  \left\{
		\begin{alignedat}{1}
			&\forall i\in\mc I_\mc P,\ c_i\leq p_i^\intercal A_K\xi_i \\
			&\forall i\in\mc I_\mc P,\ \forall j\in\mc I_\mc P,\ p_j^\intercal \xi_i\leq c_i+d^*_i
		\end{alignedat}\right. 
	\end{alignedat}
\end{equation}
with $c:=(c_1,...,c_r)^\intercal$ and $\xi:=\{\xi_i\in\Rn{n}\}_{i=1}^r$. In view of Proposition~\ref{prop:03}, $\mc R^\epsilon(q^*)$ is a probabilistic positively invariant set for the system~\eqref{eq:02:disturbanceSys} with probability $1-\epsilon$. In~\cite{Trodden2016}, both $c^*$ and $d^*$ are computed simultaneously through a single optimization problem. However, as noted in Remark 4 of~\cite{Trodden2016}, $d^*$ can also be pre-computed. In addition, a recent study~\cite{rakovic2024} has proposed a systematic way for designing the collection of norm vectors $\{p_i : \forall i \in \mc I_\mc P\}$.

\section{ROBUSTIFYING MPC IN PROBABILISTIC WAYS}
\label{sec:04}
Within our setting, given $N\in\N_+$, for any nominal state $z_k\in\Rn{n}$ at the sampling time $k\in\N$, the proposed MPC solves an OCP of the form:
\begin{equation}\label{eq:04:ocp}
	\begin{alignedat}{2}
		V^\star_N(z_k)=
		&\min_{\mbf{z}_{k},\, \mbf{v}_{k}} &&\  \sum^{N-1}_{t=0}\ell(z_{t|k}, v_{t|k}) +V_f(z_{N|k}) \\
		&\ \ \, \text{s.t.} &&\  \left\{
		\begin{alignedat}{1}
			&\forall t\in\mb I_0^{N-1},\ z_{t+1|k} = Az_{t|k}+Bv_{t|k},\\
			&\forall t\in\mb I_1^{N-1},\ z_{t|k}\in\mc Z,\\
			&\forall t\in\mb I_0^{N-1},\ v_{t|k}\in\mc V,\\
			&z_{N|k}\in\mc Z_f,	\ z_{0|k}=z_k,	
		\end{alignedat}\right. 
	\end{alignedat}
\end{equation}
with $\mbf{z}_{k}:=(z_{0|k}^\intercal,...,z_{N|k}^\intercal)^\intercal$ and $\mbf{v}_{k}:=(v_{0|k}^\intercal,...,v_{N-1|k}^\intercal)^\intercal$. The subscript $t|k$ denotes the predictions of the nominal state and control $t$-steps ahead of the sampling time $k$. The stage and terminal costs are given by
\begin{equation}
	\ell(z,v) = z^\intercal Q z + v^\intercal R v \quad\text{and}\quad 
	V_f(z) = z^\intercal P z \;,
\end{equation}
for all $z\in\Rn{n}$ and all $v\in\Rn{m}$, with $P,Q\in\mathbb{S}^{n}_{++}$ and $R\in\mathbb{S}^{m}_{++}$. The stage constraint sets on the nominal state and input are given by
\begin{equation}
	\mc Z:=\mc X\ominus\mc R^{\epsilon_x}(q^*)\quad \text{and}\quad \mc V:=\mc U\ominus K\mc R^{\epsilon_u}(q^*).
\end{equation}

\begin{assumption}\label{ass:04}
	The nominal state constraint set $\mc Z$ is a nonempty convex polyhedron that contains the origin in its interior, and the nominal input constraint set $\mc V$ is a nonempty convex polytope that contains the origin in its interior. 
\end{assumption}

This assumption implies that $\mc R^{\epsilon_x}(q^*)\subseteq\mathrm{interior}(\mc X)$ and $K\mc R^{\epsilon_u}(q^*)\subseteq\mathrm{interior}(\mc U)$, i.e.,
\begin{align*}
	\forall i\in\mc I_\mc X, \quad &\support{\mc R^{\epsilon_x}(q^*)}{f_i}<1\text{ and}\\
	\forall i\in\mc I_\mc U, \quad &\support{\mc R^{\epsilon_u}(q^*)}{K^\intercal g_i}<1.
\end{align*}
Therefore, the sets $\mc Z$ and $\mc V$ are non-empty and efficiently computed as follows:
\begin{align*}
	\mc Z &:=\left\{z\, :\, \forall i\in\mc I_{\mc X},\ f_i^\intercal z\leq 1-\support{\mc R^{\epsilon_x}(q^*)}{f_i}  \right\},\\
	\mc V &:=\left\{v\, :\, \forall i\in\mc I_{\mc U},\ g_i^\intercal v\leq 1-\support{\mc R^{\epsilon_u}(q^*)}{K^\intercal g_i}   \right\}.
\end{align*}

The nominal terminal constraint set $\mc Z_f$ and cost $V_f$ satisfy the following natural conditions~\cite{rawlings:mayne:2009}.
\begin{assumption}\label{ass:05} $\,$
\begin{enumerate}[$\,\,\ $i.]
    \item Let $K_f\in\Rnm{n}{m}$ be such that $A_{K_f}:=A+BK_f$ is strictly stabilizable. The set $\mc Z_f$ is a convex polyhedron, and it is the maximal positively invariant set for the system $z^+=A_{K_f}z$ and constraints $(z, K_fz)\in\mc Z\times \mc V$.
    \item It holds that, for all $z\in \mc Z_f$, 
    \begin{align}
        V_f((A+BK_f)z) + \ell(z,K_fz)  \leq 	V_f(z).
    \end{align}
\end{enumerate}
\end{assumption}

We define the feasible set for $N\in\N_+$ by 
\[
\mb X_N:=\left\{z_k\in\Rn{n}: z_k\text{ is such that~\eqref{eq:04:ocp} is feasible}\right\}.
\]
In our setting, the optimization problem~\eqref{eq:04:ocp} is a convex quadratic programming problem, feasible for all $z_k\in\mb X_N$. The problem involves $N(n+m) + n$ decision variables, with $Nn + n$ affine equality constraints and $N(n_\mc X + n_\mc U) + n_{\mc Z_f}$ affine inequality constraints, where $n_{\mc Z_f}$ denotes the number of the
inequality representations of the set $\mc Z_f$. We denote the parametric solution map of~\eqref{eq:04:ocp} by $\mbf{z}_k^\star(z_k)$ and $\mbf{v}_k^\star(z_k)$. Finally, we denote the MPC feedback law on the nominal system~\eqref{eq:02:nominalSys} and the real system~\eqref{eq:02:Sys} by
\begin{align*}
\bar{\kappa}_N(z_k) &= v_{0|k}^\star(z_0)\ \text{and} \\
\kappa_N(x_k) &=\bar{\kappa}_N(z_k)+K(x_k-z_k),
\end{align*}  
respectively, such that the closed-loop controlled systems are given by 
\begin{align}
	\forall k\in \N,\quad &z_{k+1} = A z_k +B\bar{\kappa}_N(z_k) \text{ and}\label{eq:04:nominalcontrolled}\\
	\forall k\in \N,\quad &x_{k+1} = Ax_k+B\kappa_N(x_k)+w_k,\ w_k\sim\omega.\label{eq:04:realcontrolled}
\end{align}  

\begin{algorithm}[htbp!]
\caption{Robustifing MPC of chance constrained linear systems}
\textbf{Initialization}: Given an initial system state $x_0$, choose an initial nominal state $z_0$ such that $x_0-z_0\in\mc R^{\epsilon_x}(q^*)$.

\textbf{for all} $k=0\to\infty$, \textbf{do}
\begin{enumerate}
    \item Solve the OCP~\eqref{eq:04:ocp} to obtain $\mbf{z}_k^\star(z_k)$ and $\mbf{v}_k^\star(z_k)$.
    \item Apply the $\bar{\kappa}_N(z_k)$ and $\kappa_N(x_k)$ to \eqref{eq:04:nominalcontrolled} and ~\eqref{eq:04:realcontrolled}, respectively. 
    \item Measure the real state $x_{k+1}$ and nominal state $z_{k+1}$ from \eqref{eq:04:nominalcontrolled} and \eqref{eq:04:realcontrolled}, respectively.
    \item Set $k\leftarrow k+1$, and go
    to Step 1.
\end{enumerate}
\end{algorithm}

The proposed MPC scheme is briefly summarized in the Algorithm 1. By construction, we have that $z_k\in\mc Z$ for all $k\in\N_+$ and $v_k\in\mc V$ for all $k\in\N$.
Thus, the condition $x_0-z_0\in\mc R^{\epsilon_x}(q^*)$, along with the state and input decompositions given in~\eqref{eq:02:controldecomposition}--\eqref{eq:02:statedecomposition}, ensures that chances constraints~\eqref{eq:02:PxkX} and~\eqref{eq:02:PukU} are satisfied.
\begin{theorem}
    Suppose Assumptions~\ref{ass:01}--\ref{ass:05} hold. Denote  $\rho_\infty$ as the stationary process of the system $s_{k+1} = A_Ks_k + w_k$, with $w_k\sim\omega$. For all $z_0\in\mb X_N$, the MPC optimization problem~\eqref{eq:04:ocp} is recursively feasible. The MPC-controlled nominal system~\eqref{eq:04:nominalcontrolled} is asymptotically stable to the origin. The real state $x_k$ of the MPC controlled system~\eqref{eq:04:realcontrolled} converges to the stationary process $\rho_\infty$  as $k\to\infty$.
\end{theorem}
\pf The recursive feasibility of~\eqref{eq:04:ocp} is guaranteed by the invariant property of the terminal set $\mc Z_f$ enforced in Assumption~\ref{ass:05}--i, as shown in~\cite{Mayne2000,rawlings:mayne:2009}.  The origin
is asymptotic stable for the dynamics~\eqref{eq:04:nominalcontrolled}, with the domain of attraction equaling to $\mb X_N$. This can
be verified by showing that the value function $V_N^0$ is a Lyapunov function for the dynamics~\eqref{eq:04:nominalcontrolled} based on Assumption~\ref{ass:05}~\cite{Mayne2000,rawlings:mayne:2009}. Because $A_K$ is strictly stabilizable, $s_k$ converges to a stationary process $\rho_\infty$ as $k\to\infty$. Given that $x_k=z_k+s_k$ and $u_k=v_k+Ks_k$, and noting that $z_k\to0$  and $v_k\to0$ as $k\to\infty$, it follows that $x_k$ converges to $\rho_\infty$. \qed

\section{NUMERICAL CASE STUDY}
\label{sec:05}
We consider a DC-DC converter model that has previously been adopted in~\cite{Cannon2011Tac,Lorenzen2017}, and it is of the form
\[
\forall k\in \N,\ x_{k+1} = \begin{bmatrix}
	1  &0.075\\
-1.43  &0.996
\end{bmatrix}x_k+\begin{bmatrix}
4.798\\
0.115
\end{bmatrix}u_k+w_k,
\]
where the mean vector and covariance matrix of the stochastic disturbance $w_k\sim\omega$ are given by 
\begin{equation}\label{eq:05:omega}
	\mu_\omega = \begin{bmatrix}
		0.005\\
		0.005
	\end{bmatrix}\quad \text{and}\quad\Sigma_\omega=10^{-4}\times\begin{bmatrix}
		1 &0\\
		0 &1
	\end{bmatrix}.
\end{equation}

The state and input constraint sets are given by
\begin{align*}
	\mc X&:=\left\{x\in\Rn{2}\ \middle|\ -2\leq[1\ 0] x \leq 2,\ -3\leq[0\ 1] x \leq 3 \right\}\\
	\mc U&:=\left\{u\in\R\ \middle|\ -0.4\leq u\leq 0.4 \right\}.
\end{align*}
The weighting matrices for the stage and terminal costs are given by
\[
Q=\begin{bmatrix}
	1 &0\\
	0 &10
\end{bmatrix}, \ R=1  \text{ and }\ P=\begin{bmatrix}
1.9074 &-5.0562\\
-5.0562 &39.5448
\end{bmatrix}.
\]
For the sake of simplicity, the control feedback and local terminal feedback matrices are specified by
\[
K=K_f=[-0.2858\ \, 0.4910],
\]  
such that the matrices $A_K$ and $A_{K_f}$ are strictly stable. Note that $K$ and $K_f$ are not required to be the same.
The terminal weighting matrix $P$ and the feedback matrix $K_f$ are derived from the solution to the infinite-horizon unconstrained optimal control problem for $(A,B,Q,R)$. The modeling parameters in chance constraints~\eqref{eq:02:PxkX}-\eqref{eq:02:PukU} are specified by
\[
\epsilon_x=\epsilon_u=0.2.
\] 
To compute a polytopic probabilistic positively invariant set given by~\eqref{eq:04:Rq}, we define the normal vectors $p_i\in\Rn{2}$ as 
\[
\forall i\in\{1,...,r\},\quad p_i = \begin{bmatrix}
	\sin \left(\frac{2\pi(i-1)}{r}\right) &\cos \left(\frac{2\pi(i-1)}{r}\right)
\end{bmatrix}^\intercal, 
\]
with $r=66$. Then, $q^*=c^*+d^*$ is computed from the solutions to~\eqref{eq:04:dstar} and \eqref{eq:04:cstar}.

\begin{figure}[htbp!]
	\centering
	\begin{overpic}[trim={1.8cm 0.5cm 1cm 1cm},clip,width=0.45\textwidth]{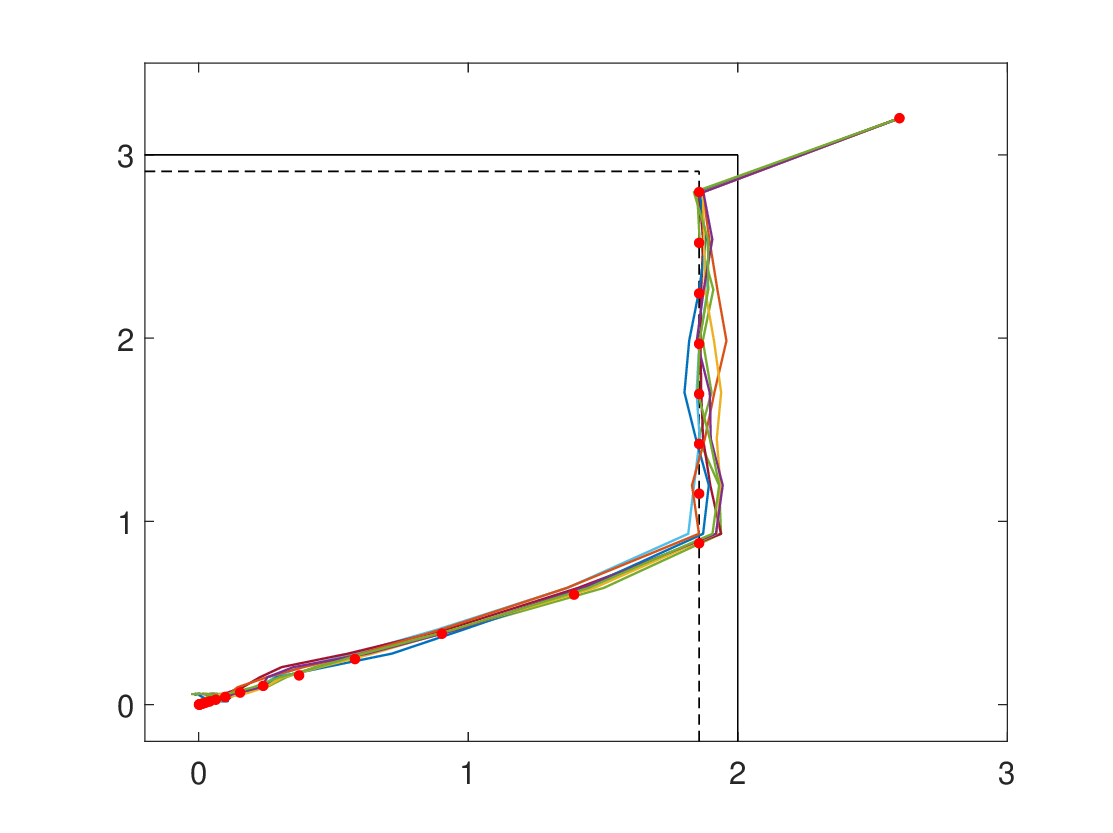}
		\put(90.5,7.2){$x^1$} \put(5.5,73){$x^2$} \put(86,71){$x_0$}
	\end{overpic}
	\caption{The closed-loop real state trajectories of~\eqref{eq:04:realcontrolled} with $8$ sampled disturbance sequences (coloured lines) and the closed-loop nominal state trajectories of~\eqref{eq:04:nominalcontrolled} (red dots). The dashed and solid lines denote partial borders of $\mc Z$ and  $\mc X$, respectively. } 
	\label{fig:01}
\end{figure}

We choose a prediction horizon $N=10$ and start the control process from an initial state $x_0 = [2.6\ \, 3.2]^\intercal$. For convenience, we set $z_0=x_0$. Fig. 1 visualizes the closed-loop real state trajectories and nominal trajectories with some disturbance sequences $\{w_k\}_{k=0}^{25}$ sampled from the multivariate normal distribution whose mean and covariance are given in~\eqref{eq:05:omega}. By simulating the closed-loop system with $10^4$ different realizations of the disturbance sequence $\{w_k\}_{k=0}^{25}$, we observe that the average state constraint violation in the first nine steps is $2\%$, while the input constraints were not violated. These results are rather conservative than the pre-defined violating probabilities $\epsilon_x=\epsilon_u=0.2$. This discrepancy is due to the fact that the computed probabilistic positively invariant sets are conservative since only mean and covariance information are used in computing the confidence regions. Fig. 2 shows the nominal state trajectory of $z_k$ given by~\eqref{eq:04:nominalcontrolled} and  the nominal input trajectory of $\bar{\kappa}_N(z_k)$, illustrating that the controlled nominal system asymptotically convergences to the origin. Consequently, the real state $x_k$ of~\eqref{eq:04:realcontrolled} converges to the stationary process, and the set $\mc R^\epsilon_x(q^*)$ is confidence region with probability $1-\epsilon_x$ of this stationary process.

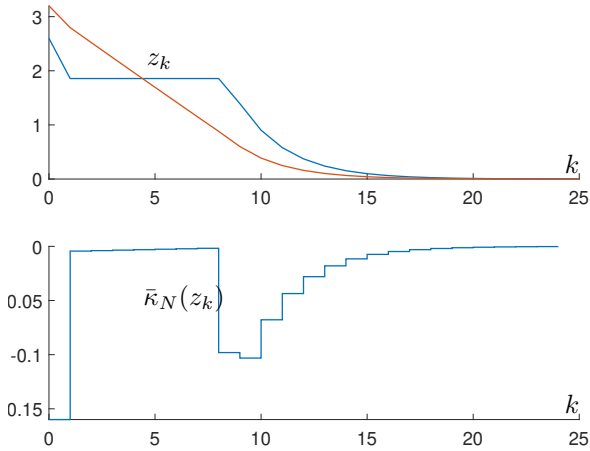
\begin{figure}[htbp!]
	\centering
	\begin{overpic}[trim={1.3cm 0.5cm 0.9cm 0.8cm},clip,width=0.45\textwidth]{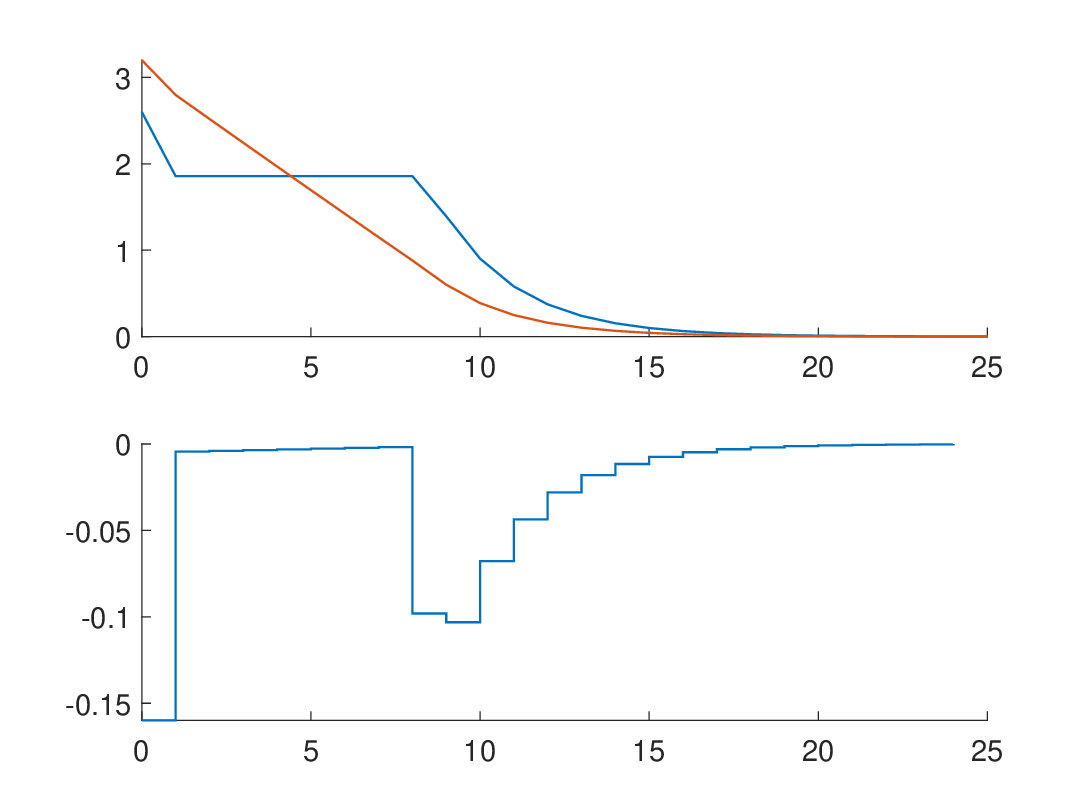}
		\put(92.5,7.2){$k$} \put(92.5,47){$k$} 
		\put(23,65){$z_k$} 
		\put(22.5,25){$\bar{\kappa}_N(z_k)$}
	\end{overpic}
	\caption{The closed-loop nominal state $z_k$ (upper part) and the closed-loop nominal input $\bar{\kappa}_N(z_k)$ (lower part).} 
	\label{fig:02}
\end{figure}

\section{CONCLUSIONS}
\label{sec:06}
This paper has introduced a model predictive controller that provides sufficient robustness to reject possibly unbounded stochastic disturbances with chance constraints on the system state and input being guaranteed. Compared to conventional MPC algorithms, the proposed method only needs minor offline computational efforts to compute a probabilistic positively invariant set, which can be easily computed by solving a simple linear programming problem. Future work will aim to reduce conservativeness, adapt the method for nonlinear systems, and extend the certainty-equivalent cost function to expected ones.


%
%

\bibliographystyle{IEEEtran}
\bibliography{references.bib}

\end{document}